\documentclass[preprint,12pt]{elsarticle}

\usepackage{amssymb}

\usepackage{amsmath,amsfonts,amscd}
\usepackage[all,cmtip]{xy}
\usepackage{eucal}
\AtBeginDocument{%
\newwrite\XFilenames%
\immediate\openout\XFilenames\jobname.x-fls%
\immediate\write\XFilenames{#-BEGIN-Filenames}%
}
\AtEndDocument{%
\immediate\write\XFilenames{#-END-Filenames}%
\immediate\closeout\XFilenames%
}

\def\hideDETAIL#1{}

\numberwithin{equation}{section}

\newtheorem{theorem}{Theorem}[section]
\newtheorem{proposition}[theorem]{Proposition}
\newtheorem{lemma}[theorem]{Lemma}
\newtheorem{corollary}[theorem]{Corollary}

\newproof{pf}{Proof}
\newdefinition{remark}{Remark}[section]
\newdefinition{definition}{Definition}[section]

\def\HilH{{\mathcal H}}

\def\AA{{\mathcal A}}
\def\BB{{\mathcal B}}
\def\DD{{\mathcal D}}
\def\CC{{\mathbb C}}
\def\EE{{\mathcal E}}
\def\FF{{F}}

\def\KK{{\mathcal K}}

\def\RR{{\mathbb R}}

\def\XX{{\mathcal X}}

\def\YY{{\mathcal Y}}

\def\la{\langle}
\def\ra{\rangle}
\def\DOM{{\mathcal D}}

\def\NULL{\mathcal N}
\def\RANGE{\mathcal R}
\def\IM{\operatorname{Im}}

\def\SPEC{\sigma}
\def\DISC{\operatorname{disc}}
\def\RESOLVENT{\rho}
\def\CLOSURE#1{\overline{#1}}
\def\ESS{\operatorname{ess}}

\def\AC{{\operatorname{ac}}}

\def\QED{\hfill$\square$\medskip}

\def\DeferNotation#1{}

\journal{}

\begin{document}

\begin{frontmatter}


\author[tyn]{Take-Yuki ~Nagao}
\ead{nagao-takeyuki@aiit.ac.jp}
\address[tyn]{Advanced Institute of Industrial Technology,
1-10-40 Higashi-Ohi, Shinagawa-City, Tokyo, JAPAN 140-0011
}

\title{
Generalized Fourier representation
of the absolutely continuous part of
a selfadjoint operator
}




\begin{abstract}
We formulate and prove the existence and uniqueness of
the generalized Fourier transform associated with
the absolutely continuous part of an arbitrary
selfadjoint operator on a separable Hilbert space.
To this end we develop a novel method to
decompose an absolutely continuous operator
into a variable fiber direct integral of 
selfadjoint operators.
\end{abstract}


\end{frontmatter}


  \section{Introduction}
%
%
%
The most fundamental result of the spectral theory of 
selfadjoint operators
is the spectral theorem which states
$$
H = \int_\RR \lambda d E_\lambda
$$
where $H$ is a selfadjoint operator on some Hilbert space $\HilH$
and $\{E_\lambda\}_{\lambda\in\RR}$ is a spectral family associated with $H$.
%
In this paper we establish a new spectral representation
theorem that connects
the spectral family and the generalized eigenfunctions.

By generalized eigenfunction we mean a solution $u$
to the formal eigenequation
\begin{equation}
\label{_label_luyvnsdzafgqbcbu}
(H-\lambda)u = 0,\quad u \in \XX
\end{equation}
where $\lambda \in \RR$ is a spectral parameter.
The function space $\XX$ must be chosen
suitably, depending on the operator $H$ of interest.
Basically we choose $\XX$ 
in such a way that the solutions to
the problem (\ref{_label_luyvnsdzafgqbcbu})
carry a complete set of information about the operator $H$
at the energy $\lambda$.
We can take $\XX = \HilH$ when $\lambda \notin \SPEC(H)$
or $\lambda \in \SPEC_{\DISC}(H)$.
But we need to have $\XX$ larger than $\HilH$
if $\lambda \in \SPEC_{\ESS}(H)$
in order to capture the generalized eigenfunctions
corresponding to the continuum.

The goal of this paper is show the existence of
a family of Hilbert spaces
$\{\XX(\lambda)\}_{\lambda \in \RR}$
such that
there is a selfadjoint operator $H_\lambda$ on
each $\XX(\lambda)$
with the same spectrum as $H$
and the collection of the solutions to the problem
\begin{equation}
\label{_label_mytjrqcvphcloozl}
(H_\lambda - \lambda)u = 0,
\quad u \in \DOM(H_\lambda)
\end{equation}
for $\lambda$ from a set of full Lebesgue
measure determines the absolutely continuous part
of $H$. We shall show that the spectral density $E_\lambda'$
is a pullback of the orthogonal projection
$Q(\lambda)$
onto the null space of $H_\lambda - \lambda$
in $\XX(\lambda)$ in the sense that
\begin{equation}
\label{_label_sfovdosnpgtinnva}
\la E'_\lambda f | f \ra
= \|Q(\lambda)Y(\lambda)f\|^2,
\quad
f \in \DOM(Y(\lambda)),
\end{equation}
where 
$Y(\lambda) : \HilH \longrightarrow \XX(\lambda)$
is a possibly non-closable injective linear map.
This means that the generalized eigenvalue problem
(\ref{_label_mytjrqcvphcloozl})
determines the absolutely continuous part of $H$.

We now briefly compare the idea of this paper with the known results.
In the case of partial differential opeartors,
it is common to use weighted $L^2$ spaces \cite{MR53:1053}
and Besov spaces
\cite{springerlink:10.1007/BF02786703%
,HoermanderII%
,springerlink:10.1007/PL00005579%
} for the function space $\XX$ in 
(\ref{_label_luyvnsdzafgqbcbu}).
Well known technique is the following.
Given an absolutely continuous operator
$H$ on some function space, we look for
a function space $\YY$ such that
we have
\begin{equation*}
\label{_label_ijnehchebidmzvdd}
|\la (H-z)^{-1} f | g \ra|
\leq C \|f\|_{\YY}\|g\|_{\YY},
\quad
f, g \in \YY
\end{equation*}
with a uniform constant $C > 0$.
The choice of $\YY$ depends on the problem in general.
The point is that we have
\begin{equation*}
\label{_label_zgrxwbjvmrhtaesg}
\IM \la (H-z)^{-1}f | f \ra
\leq C\|f\|_{\YY}^2,
\quad
f \in \YY,
\quad \IM z > 0
\end{equation*}
so that there is a bounded operator
$\delta(H-\lambda) : \YY \longrightarrow \YY^*$
such that
\begin{equation}
\label{_label_okahzmqilbfldjuu}
\lim_{\varepsilon\downarrow 0}
\IM \la (H-\lambda -i\varepsilon)^{-1}f | f \ra
= \pi \la \delta(H-\lambda) f | f \ra,
\quad
f \in \YY
\end{equation}
for all $\lambda \in \RR$.
The operator
$\delta(H-\lambda)$ is identical to the
spectral density $E_\lambda'$ realized as an
opeartor.
One usually fix a function space $\DD$ and
argues that any solution $u \in \YY^*$ of
the problem
\begin{equation}
\label{_label_ggkclwktsjmuitjz}
\la u | (H-\lambda)f \ra = 0,
\quad
f \in \DD
\end{equation}
can be written as
$
u = \delta(H-\lambda) \varphi
$
for some $\varphi \in \YY$ and that we have
\begin{equation}
\label{_label_ecoxyejzoppzcjlf}
\la Hf | f \ra
= \int_\RR \lambda \la \delta (H-\lambda) f | f\ra d\lambda.
\end{equation}
The identity
(\ref{_label_ecoxyejzoppzcjlf})
guarantees that the collection of solutions to
the problem
(\ref{_label_ggkclwktsjmuitjz})
for all $\lambda \in \RR$ is complete in the sense
that we can recover the original operator $H$
from this collection.

The main idea of this paper is to realize the left hand
side of
(\ref{_label_sfovdosnpgtinnva}),
which roughly corresponds to
(\ref{_label_okahzmqilbfldjuu}),
as a quadratic from and choose the completion of its graph
as $\YY$.  Although the quadratic form is nonnegative,
it is not always closable.  We shall establish
a new representation theorem of nonnegative quadratic
forms to overcome this difficulty.
Another problem is how to give a rigorous meaning
to the problem
(\ref{_label_ggkclwktsjmuitjz}).
To solve this,
we propose the notion of continuation of a selfadjoint 
operator to some other function spaces.

The author expects that the results of this paper is useful
in scattering theory, since one can prove
the existence of the generalized Fourier transform,
which is a key tool for the theory,
without depending on the limiting absorption
principle nor on the Mourre inequality.

In this paper
we focus on the absolutely continuous part of $H$ and
the results on the case where $H$ has a non-trivial
singular continuous part will appear elsewhere.

We close this section by collecting the notations used throughout
the paper.
A quadratic form $q$ is linear in the first argument
and anti-linear in the second. We adopt the same convention
about inner products.
We write $q(f)=q(f,f)$.
By
$\RESOLVENT(T)$,
$\DOM(T)$,
$\RANGE(T)$
and
$\NULL(T)$
we mean
the resolvent set, domain, range and null space of
a linear opeartor $T$.
If $T$ is closable we write $\CLOSURE{T}$ for its closure
unless otherwise noted.
By $P_\AC(H)$ we denote the spectral projection
onto the absolutely continuous subspace of a selfadjoint operator $H$.
The set of infinitely differentiable complex valued
functions $f$ on $\RR$ such that $f$ and all its derivatives
tend to zero at infinity is denoted by $C_0^\infty(\RR)$.
We sometimes write $0$ for a zero space.  The identity operators
of various function spaces are denoted by $1$.
$T\restriction \DD$ is the operator $T$ with
domain restricted to a subspace $\DD$ of $\DOM(T)$.
  \section{Continuation}
In this section we develop a method to continue a selfadjoint operator
$H$
initially defined on a Hilbert space $\HilH$ to another Banach space
$\XX$ that is close to $\HilH$ in a certain sense.
We will use the results of this section to formulate the
generalized eigenspace of a selfadjoint operator.
Throughout this section we assume that $\XX$ is a Banach space
obtained by the completion of a dense subspace $\DD$ of $\HilH$
with respect to some norm on $\DD$.

The goal of this section is to formulate a condition
that guarantees the existence of
a closed operator
$H_\XX : \XX \longrightarrow \XX$,
which we will call a continuation of $H$ to $\XX$.
The operator $H_\XX$ will be chosen so that
it coincides with $H$ on a common core of
$H_\XX$ and $H$.

We start with the observation that
$\XX$ is isometrically isomorphic to the completion of $\DD$
by some norm
if and only if there is a densely defined injective linear map
$
Y : \HilH \longrightarrow \XX
$
with dense range
such that
$
\DOM(Y) = \DD.
$
In this paper we call such $Y$ a completion operator.
Note that a completion operator need not be closable.
Just like in the case of closed operators,
we say that a subspace $\DD$ of $\XX$ is
a core of $Y$ if for any $f \in \DOM(Y)$ there is a
sequence $f_j \in \DD$ such that $f_j \rightarrow f$
and $Yf_j \rightarrow Yf$.
In this case the restriction $Y\restriction \DD$
is a completion operator from $\HilH$ to $\XX$.

Given a completion operator $Y$,
we someteimes identify $\DOM(Y)$ with $\RANGE(Y)$
through the correspondence
$\DOM(Y)\ni f \mapsto Yf \in \RANGE(Y)$
in order to simplify notations.
We also remark that since $\XX$ is a completion,
it cannot always be identified with a subspace
of $\HilH$.  As we will see later,
however, it is important to consider such
a case when dealing with a non-closable form.

Our strategy
is that we first continue the resolvent $(H-z)^{-1}$
to $\XX$ and then define the continuation $H_\XX$ using
the inverse of the continued resolvent.  We introduce
the following
\begin{definition}
\label{_label_zguvvrankpkeptbb}
Let $Y : \HilH \longrightarrow \XX$ be a 
completion operator and $z \in \RESOLVENT(H)$.
We say that $H$ is compatible with $Y$ at $z$
if there is a subspace $\DD \subset \DOM(Y)$
with the following properties:
\begin{enumerate}[i.]
\item $(H-z)^{-1}\DD$ is a core of $Y$;
\item there exists a constant $C_z > 0$ depending of $z$ such that
$$
\|Y(H-z)^{-1}f\|_\XX
\leq C_z \|Yf\|_\XX
$$
for all $f \in \DD$;
\item $u \in \XX$, $f_j\in\DD$, $Yf_j \rightarrow u$ in $\XX$,
and $Y(H-z)^{-1}f_j \rightarrow 0$ in $\XX$ imply $u = 0$.
\end{enumerate}
Given a subset $\Omega \subset \RESOLVENT(H)$,
we say that $H$ is compatible with $Y$ in $\Omega$
if $H$ is compatible with $\XX$ through $Y$
for all $z \in \Omega$.
\end{definition}

\begin{theorem}
\label{_label_btxjzynpxykhuqck}
Let $\Omega$ be a non-empty subset of $\RESOLVENT(H)$.
Suppose that $H$ is compatible with a completion operator
$Y:\HilH \longrightarrow\XX$ in $\Omega$.
Let $\DD$ be as in Definition \ref{_label_zguvvrankpkeptbb}.
Then there exists a densely defined closed operator
$H_\XX : \XX \longrightarrow \XX$ such that
$\Omega \subset \RESOLVENT(H_\XX)$ and that
\begin{equation}
\label{_label_dzcvxwxsprmxhxkb}
(H_\XX -z)^{-1}Yf = Y(H-z)^{-1}f
\end{equation}
for all $f \in \DD$ and  $z \in \Omega$.
The operator $H_\XX$ is independent of the choice of
$\DD$ and uniquely determined by $H$ and $Y$.
\end{theorem}
\begin{pf}
There exists a unique bounded
operator
$B(z)$ on $\XX$ such that
\begin{equation}
\label{_label_ottclibwtjhkrarb}
B(z)Yf = Y(H-z)^{-1}f
\end{equation}
for all $f \in \DD$.
The first resolvent identity of $(H-z)^{-1}$
yields
\begin{equation}
\label{_label_ygsmvnwcclziaxdw}
B(z) - B(w) = (z-w)B(z)B(w)
\end{equation}
for all $z, w \in \Omega$ and it follows that
$\RANGE(B(z))$ is constant for all $z \in \Omega$.
Remark that
$B(z)$ is injective and has a dense range.
We define
the operator $H_\XX : \XX \longrightarrow \XX$
with domain $\DOM(H_\XX) = \RANGE(B(z))$
by
\begin{equation}
\label{_label_vypgyobhxpjsiqtp}
H_\XX B(z) f = f + zB(z)f,\quad
f \in \XX.
\end{equation}
Since $B(z)$ is bounded and injective, $H_\XX$ is closed.
Moreover,
$H_\XX$ is independent of the choice of $z \in \Omega$
due to (\ref{_label_ygsmvnwcclziaxdw}).
The identity (\ref{_label_vypgyobhxpjsiqtp}) implies that
$H_\XX-z$ is bijective and satisfies $(H_\XX-z)^{-1} = B(z)$
for all $z \in \Omega$.
The uniqueness of $H_\XX$ follows from the uniqueness of
$\BB(z)$ in (\ref{_label_ottclibwtjhkrarb}).
This completes the proof.
\QED
\end{pf}

When we identify $f \in \DOM(Y)$ and $Yf \in \RANGE(Y)$,
(\ref{_label_dzcvxwxsprmxhxkb}) means that
$H$ and $H_\XX$ share 
$\DD_0 = (H-z)^{-1}\DD = (H_\XX-z)^{-1}\DD$
as a common core and that
we have
$H_\XX u = H u$ for all $u \in \DD_0$.
It is thus reasonable to regard $H_\XX$ as an
extension of $H$ to the Banach space $\XX$.
In order to avoid confusion with the usual notion of
extension of operators,
let us refer to the operator $H_\XX$ as
a continuation of $H$ to $\XX$ induced by $Y$
or simply a continuation of $H$.
In general, $H_\XX$ depends on the choice of $Y$
and not uniquely determined by $H$ and $\XX$.

\begin{corollary}
\label{_label_qahdlsdgvfgxtoqj}
Suppose that $\XX$ is a Hilbert space and that
$H$ is compatible with a completion operator
$Y : \HilH \longrightarrow \XX$ in $\{i, -i\}$.
Suppose further that $H_\XX$ is symmetric.
Then $H_\XX$ is selfadjoint.
\end{corollary}
  \section{Pullback representation of a nonnegative form}
In this section we shall develop our main tool
to deal with nonnegative quadratic forms that might not be
closable.  Our goal here is to formulate and prove
a representation theorem of nonnegative quadratic forms.
Recall that a nonnegative quadratic form $q$
on a Hilbert space $\HilH$
is closable if and only if $\DOM(q)$ is complete
under the norm $\|\cdot\|_q$ where
$
\|f\|_q^2 = \|f\|^2 + q(f).
$
So it is reasonable to embed the domain $\DOM(q)$
into its completion by the above norm.
We thus begin with the following
\begin{definition}
Suppose that $\HilH, \XX$ and $\EE$ are Hilbert spaces.
A pullback triple $\tau=(J, F, Y)$ is a triplet of operators
\begin{equation}
\label{_label_ugkdvuhlcbthklim}
J : \XX \longrightarrow \HilH,
\quad
\FF : \XX \longrightarrow \EE,
\quad
Y: \HilH \longrightarrow \XX
\end{equation}
with the following properties:
\begin{enumerate}[i.]
\item $J, F$ are bounded with dense range;
\item $Y$ is a completion operator;
\item for any $u \in \XX$ we have
$
\|u\|_{\XX}^2 = \|Ju\|^2 + \|Fu\|_{\EE}^2;
$
\item
\label{_label_eigyitvlgjudsfdc}
$JY \subset 1$.
\end{enumerate}
We often say that $\tau$ is a pullback triple on $\HilH$
in order to specify the function space $\HilH$ explicitly.
We call $\DOM(Y)$ the domain of $\tau$ and denote it by $\DOM(\tau)$.
\end{definition}
Note that the condition
(\ref{_label_eigyitvlgjudsfdc})
is equivalent to the commutativity of
the diagram
$$
\xymatrix{
&\DOM(Y) \ar[ld]_{\iota} \ar[d]^{Y}
\\
\HilH &\ar[l]_{J} \XX
}
$$
where $\iota$ is the natural injection.
As is done in depicting the operator $Y$ of the above diagram,
we put the domain of a linear operator
in front of the tail of the arrow indicating the operator.
We prefer this
rule in order to simplify notations.

\begin{lemma}
\label{_label_pwuxxmcauqsiqpoj}
Let $(J, F, Y)$ be a pullback triple.
Then
\begin{enumerate}[i.]
\item $Y^{-1}$ is closable and $J = \CLOSURE{Y^{-1}}$;
\item
$Y$ is closable if and only if $J$ is injective
and in this case we have $\CLOSURE{Y} = J^{-1}$;
\item
$Y$ is bounded if and only if $J$ is an isomorphism.
\end{enumerate}
\end{lemma}
\begin{pf}
(i).  Since $JYf = f$ for any $f \in \DOM(Y)$,
we have $\|Yf\|^2_{\XX} = \|f\|^2 + \|FYf\|_{\EE}^2\geq \|f\|^2$ for
any $f \in \DOM(Y)$.  It follows that
$Y^{-1}$ is closable and $J = \CLOSURE{Y^{-1}}$.
(ii).
If $f_j \in \DOM(Y)$ and $Yf_j \rightarrow u$ in $\XX$,
then we have $f_j \rightarrow Ju$ in $\HilH$.
It follows that $Y$ is closable if and only if $J$ is injective.
It is easy to see that $\DOM(\CLOSURE{Y}) = \RANGE(J)$
and that $\CLOSURE{Y}Ju = u$ for all $u \in \XX$.
(iii).
Note that $\|Yf\| \leq C \|f\|$ for all $f \in \DOM(Y)$
is equivalent to $\|u\| \leq C \|Ju\|$ for all $u \in \XX$.
Since $J$ is bounded and has a dense range, the statement follows.
\QED
\end{pf}

\begin{lemma}
\label{_label_khnpkbumclnccrmi}
Let
$$
J_j : \XX_j \longrightarrow \HilH_j,
\quad
\FF_j : \XX_j \longrightarrow \EE_j,
\quad
Y_j: \HilH_j \longrightarrow \XX_j,
\quad
j = 1,2
$$
be pullback triples.
Then the following conditions are equivalent:
\begin{enumerate}[i.]
\item $J_1$ and $J_2$ are unitarily equivalent;
\item $F_1$ and $F_2$ are unitarily equivalent;
\item there exist unitary isomorphisms
$W_\HilH, W_\XX,$ and $W_\EE$ such that
the following diagram commutes.
$$
\xymatrix{
\HilH_1 \ar[d]^{W_\HilH}
  & \XX_1 \ar[l]_{J_1} \ar[r]^{F_1} \ar[d]^{W_\XX}
  &\EE_1 \ar[d]^{W_\EE}
  \\
\HilH_2 & \ar[l]_{J_2} \XX_2 \ar[r]^{F_2} &\EE_2 \\
}
$$
\end{enumerate}
\end{lemma}

We introduce the following equilvalence relation.
\begin{definition}
We say that two pullback triples $(J_j, \FF_j, Y_j), j=1,2$
are equivalent if any of the conditions of Lemma
\ref{_label_khnpkbumclnccrmi}
is satisfied.
\end{definition}

\begin{lemma}
\label{_label_iucbdqsouvzpoyfk}
Two pullback triples $(J_j, \FF_j, Y_j), j=1,2$ on a Hilbert space $\HilH$
are equivalent if and only if there is a unitary isomorphism
$W$ such that $(J_1, F_1, WY_2)$ is a pullback triple.
\end{lemma}

\begin{lemma}
\label{_label_etripfglmgjuzhth}
Let $(J, F, Y)$ be a pullback triple and $Y_0 \subset Y$.
Then $(J, F, Y_0)$ is a pullback triple if and only if
$\DOM(Y_0)$ is a core of $Y$.
In this case $(J, F, Y)$ and $(J, F, Y_0)$ are equivalent.
\end{lemma}

We say that a linear space $\DD$ is a core of a pullback triple
$(J, F, Y)$ if $\DD$ is a core of $Y$.

\begin{theorem}
\label{_label_dncrxfyquumwagzk}
Let $q$ be a densely defined nonegative quadratic form
on a Hilbert space $\HilH$.  Then there exists a pullback triple
$(J, F, Y)$ with $\DOM(Y) = \DOM(q)$ such that
the following diagram commutes.
$$
\xymatrix{
&\DOM(q) \ar[ld]_{\iota} \ar[d]^{Y} \ar[rd]
\ar[r]^{q(\cdot)}
& \RR \\
\HilH &\ar[l]_{J} \XX \ar[r]^{F}& \EE \ar[u]_{\|\cdot\|_{\EE}^2}
}
$$
The pullback triple is unique up to equivalence.
Moreover, $q$ is closable if and only if so is $Y$.
In this case $(J, F, \CLOSURE{Y})$ is a
pullback triple of $\CLOSURE{q}$ and we have
\begin{equation}
\label{_label_fuqxghaeeppppuke}
\DOM(\CLOSURE{q}) = \RANGE(J)
,\quad
\CLOSURE{q}(Ju) = \|\FF u \|_\EE^2,
\quad
u \in \XX.
\end{equation}
\end{theorem}
\begin{remark}
The commutativity of the above diagram implies
$$
q(f) = \|FYf\|_{\EE}^2,
\quad
f \in \DOM(Y).
$$
This means that any nonnegative form $q$
is a pullback of a bounded nonnegative form by
a completion operator $Y$.
\end{remark}
\begin{pf}
It is easy to see that
$$
\DOM(q)\times \DOM(q) \ni(f,g) \mapsto \la f | g \ra
 + q(f, g) \in \CC
$$
is a positive definite inner product on $\DOM(q)$.
Let $\XX$ be the completion of $\DOM(q)$ by this inner product
and $Y : \HilH \longrightarrow \XX$ be the associated
completion operator.
Since $q$ is nonnegative,
we have
$$
|q(f, g)|^2 \leq q(f) q(g) \leq \|Yf\|_\XX^2 \|Yg\|_\XX^2
$$
for any $f, g \in \DOM(q)$.  By the Riesz representation
theorem there exists a unique bounded selfadjoint operator
$\delta : \XX \longrightarrow \XX$ such that
$$
\la \delta Yf | Yg \ra_{\XX} = q(f, g),
\quad f, g \in \DOM(q).
$$
It is easy to see that the mapping
$$
\RANGE(\delta) \times \RANGE(\delta)
\ni (\delta u , \delta v) \mapsto \la \delta u | v \ra_{\XX} \in \CC
$$
is well-defined and is a positive definite inner product on $\RANGE(\delta)$.
Let $\EE$ be the completion of $\RANGE(\delta)$ with respect to this
inner product.
We define the operator $\FF$ by
$$
\FF : \XX \longrightarrow \EE,
\quad \FF u = \delta u, \quad u \in \XX.
$$
By our choice of $\EE$, $\RANGE(F)$ is dense in $\EE$
and we have
$
\|\FF u\|_\EE^2 = \la \delta u | u \ra_{\XX}
$
for any $u \in \XX$.  Note that $Y^{-1}$
is closable and has a bounded extension by
Lemma \ref{_label_pwuxxmcauqsiqpoj}.
We set $J = \CLOSURE{Y^{-1}}$.
Clearly,
\begin{equation}
\label{_label_expsuxspcvpknyal}
JYf = f,
\quad
q(f) = \|FYf\|^2_{\EE},
\quad
f \in \DOM(Y).
\end{equation}
Therefore $(J, F, Y)$ has the desired properties.
The uniqueness
follows from
(\ref{_label_expsuxspcvpknyal})
and Lemma
\ref{_label_khnpkbumclnccrmi}%
.
The statement about closability follows from Lemma
\ref{_label_pwuxxmcauqsiqpoj}%
.
The relation (\ref{_label_fuqxghaeeppppuke}) is obvious.
This completes the proof.
\QED
\end{pf}

We refer to the triple $(J, F, Y)$ as a pullback triple of
$q$ and we say that two nonnegative forms $q_1$ and $q_2$
are equivalent if the corresponding pullback triples
are equivalent.  There are cases where the domains of
$q_1$ and $q_2$ have no inclusion relation but $q_1$
and $q_2$ are equivalent.

Similarly any linear operator is a pullback
of a bounded operator.  We only state the result
and the proof is left to the reader.
\begin{theorem}
\label{_label_shjgervlestqfqkc}
Let $T : \HilH \longrightarrow \KK$ be a linear 
operator between Hilbert spaces.
Then there exists a pullback triple $(J, F, Y)$
with $\DOM(Y) = \DOM(T)$
such that the following diagram commutes.
$$
\xymatrix{
&\DOM(T) \ar[ld]_{\iota} \ar[d]^{Y} \ar[rd]^{T}
\\
\HilH &\ar[l]_{J} \XX \ar[r]^{F}& \CLOSURE{\RANGE(T)}
}
$$
The pullback triple is unique up to equivalence.
Moreover, $T$ is closable if and only if so is $Y$.
In this case $(J, F, \CLOSURE{Y})$ is a
pullback triple of $\CLOSURE{T}$ and we have
$$
\DOM(\CLOSURE{T}) = \RANGE(J),
\quad
\CLOSURE{T}Ju = Fu,
\quad u \in \XX.
$$
\end{theorem}

Fix a pullback triple
(\ref{_label_ugkdvuhlcbthklim}).
Then
$\NULL(J)$ and $\NULL(F)$ are orthogonal in $\XX$.
Moreover,
\begin{gather}
\NULL(J) = \{ f \in \XX : f = Q f\},
\quad
\label{_label_ekpkdvmkmvzchkga}
\NULL(F) = \{ f \in \XX : f = \overline{Q} f \},
\end{gather}
where $Q$ and $\overline{Q}$ are bounded operators in $\XX$
defined by
\begin{gather*}
Q = F^*F,
\quad \overline{Q} = 1 - Q = J^*J.
\end{gather*}
\begin{lemma}
\label{_label_yrubbzipptagtbxr}
Let $(J, F, Y)$ be a pullback triple.
Then the following conditions are equivalent:
\begin{enumerate}[i.]
\item $F$ is a partial isometry;
\item $\XX = \NULL(J) \oplus \NULL(F)$;
\item $JF^* = 0$;
\item $J\NULL(F)$ is dense in $\HilH$.
\end{enumerate}
\end{lemma}

A pullback triple
$(J, F, Y)$ is said to be splitting
if any of the conditions of
Lemma
\ref{_label_yrubbzipptagtbxr}
is satisfied.
We close this section by summarizing the properties of
a splitting pullback triple for later use.
\begin{proposition}
\label{_label_cqwlrlcwckueubdr}
Suppose that $(J, F, Y)$ is a splitting pullback triple.
Then $J$ and $F$ are partial isometries with
\begin{gather}
\label{_label_xrthunorxyokpopj}
\NULL(F) = \RANGE(J^*),
\quad
\RANGE(J) = \HilH,
\\
\NULL(J) = \RANGE(F^*),
\quad
\RANGE(F) = \EE.
\end{gather}
\end{proposition}
\begin{remark}
Proposition \ref{_label_cqwlrlcwckueubdr} implies
that under the condition that the pullback triple $(J, F, Y)$
is splitting, the form $q$ associated with the triple is closable if 
and only if $F$ and thus $q$ vanishes identically.
In this case we have $\EE = 0$.  Note also that
$\XX$ is unitarily isomorphic to
the splitting direct sum $\HilH \oplus \EE$.  This means that
one can measure by $\EE$ how far is $q$
from being closable.
\end{remark}
  \section{Pullback representation of spectral densities}
Suppose that $H$ is a selfadjoint operator on a separable
Hilbert space $\HilH$.
Let $\{E_\lambda\}_\lambda$ be a spectral family of $H$.
We shall derive a representation of
the absolutely continuous part of $H$ in terms of the spectral
density $E'_\lambda$.  The main diffculty is that
the spectral density is never closable unless it is zero identically.
The pullback triples
we have developed in the previous section
facilitate the understanding of
the spectral density as a quadratic form.

Since
$\|E_\lambda f\|^2$ is a non-decreasing function of $\lambda$
for any $f \in \HilH$,
it is differentiable almost everywhere
in the sense of Dini, and so is
$
\la E_\lambda f | g \ra
$
for any $f, g \in \HilH$.
It is easy to see that
\begin{equation}
\label{_label_kxhamrmbacjtwcsw}
\frac{
d\la E_\lambda f | g \ra
}{
d\lambda
}
=
\frac{
d\la E_\lambda P_{\AC}(H)f | g \ra
}{
d\lambda
}
\end{equation}
almost everywhere.  By Radon-Nikodym theorem we see that
(\ref{_label_kxhamrmbacjtwcsw})
is integrable and that
\begin{equation}
\label{_label_ftgqoznmugmufdat}
\la P_\AC(H) f | g \ra
=
\int_\RR \frac{d\la E_\lambda f | g \ra}{d\lambda}
d\lambda
\end{equation}
for any $f, g \in \HilH$.
We shall show that 
(\ref{_label_kxhamrmbacjtwcsw})
is well-defined as a densely defined nonnegative form.

\begin{lemma}
\label{_label_wbdmetvzzvzddjpo}
Let $f, g \in \HilH$.  Suppose that $\la E_\lambda f | g \ra$
is differentiable at $\lambda$.  Then
$\la E_\lambda \chi(H) f | g \ra$ is differentiable at $\lambda$
for any $\chi \in C^\infty_0(\RR)$ and we have
\begin{equation}
\label{_label_ejpglsvsbnpgnvua}
\frac{
d\la E_\lambda \chi(H)f | g \ra
}{
d\lambda
}
=
\chi(\lambda)
\frac{
d\la E_\lambda f | g \ra
}{
d\lambda
}.
\end{equation}
\end{lemma}
\begin{pf}
The statement follows from the fact that
$
[\chi(H) - \chi(\lambda)]E_H(I) / |I|
\rightarrow 0
$
strongly as $|I| \rightarrow 0$ for any
interval $I$ containing $\lambda$ as an internal point.
\QED
\end{pf}

Note that the spectral density has the
following representation by the boundary
value of the resolvent.
\begin{lemma}
Let $f \in \HilH$.  Then
\begin{equation}
\label{_label_hovikyqnniucesxg}
\lim_{\varepsilon \downarrow 0}
\frac{\varepsilon}{\pi}
\|(H-\lambda-i\varepsilon)^{-1}P_\AC(H)f\|^2
=
\frac{
d\|E_\lambda f \|^2
}{
d\lambda
}
\end{equation}
for almost every $\lambda$ with respect to Lebesgue measure.
\end{lemma}

Suppose that $\DD_0$ is a subset of a Hilbert space $\HilH$
and $\AA$
is a family of bounded linear operators on $\HilH$.
Let $\DD$ be a linear space of
all finite linear combinations of the vectors
of the form
$
A f,
f \in \DD_0,
$
where $A$ is the product of zero or finitely
many elements of $\AA$.  It is easy to see that
$
\DD_0 \subset \DD
$
and $A \DD \subset \DD$
for all $A \in \AA$.  We call $\DD$ 
the span of $\DD_0$ by $\AA$.
We are interested in the case
where $\DD_0$ is an
orthonormal basis of $\HilH$ and
$\AA = \{(H-z)^{-1}\}_{z \in \RESOLVENT(H)}$.
In this case
we say that $\DD$ is
the span of an orthonormal basis by the resolvent of $H$.



\begin{lemma}
\label{_label_ifzcfacnwzeinxea}
Let $\DD$ be the span of an orthonormal basis by the resolvent
of $H$.
Then
there exists a subset $\Lambda \subset \RR$
of full Lebesgue measure such that
$
\la E_\lambda f | g \ra
$
is differentiable at any $\lambda \in \Lambda$ whenever $f, g \in \DD$.
\end{lemma}
\begin{remark}
We remark that $\Lambda$ does not contain an eigenvalue of $H$.
\end{remark}
\begin{pf}
Let $\{\varphi_j\}$ be the orthonormal basis generating $\DD$,
which is countable by assumption.  There is a subset $\Lambda \subset \RR$
of full Lebesgue measure such that
$
\la E_\lambda \varphi_j | \varphi_k \ra
$
is differentiable at any $\lambda \in \Lambda$ for all $j, k$.
By Lemma
\ref{_label_wbdmetvzzvzddjpo}
$
\la E_\lambda A\varphi_j | B\varphi_k \ra
$
is differentiable at any $\lambda \in \Lambda$ for all $j, k$,
where $A$ and $B$ are zero or finitely many products of operators
in $\{(H-z)^{-1}\}_{z \in \RESOLVENT(H)}$.
The statement follows by linearity.
\QED
\end{pf}

\begin{lemma}
\label{_label_brlgehxetoirlenm}
Suppose that each of $\DD_1$ and $\DD_2$
is the span of an orthonormal basis by the resolvent of $H$.
Then there exists a set $\Lambda_0$ of full Lebesgue measure
with the following properties:
\begin{enumerate}[i.]
\item
$\la E_\lambda f | g \ra$ is differentiable at $\lambda$
whenever $\lambda\in \Lambda_0$ and $f, g \in \DD_1 + \DD_2$;
\item
for any $f \in \DD_2$ there exists a sequence $f_j \in \DD_1$
such that 
\begin{equation}
\label{_label_mtknmljxgqmjihhu}
\|f - f_j\|\rightarrow 0,
\quad
\frac{
d\|E_\lambda(f-f_j)\|^2
}{
d\lambda
}
\rightarrow 0
\end{equation}
as $j \rightarrow \infty$ for any $\lambda \in \Lambda_0$;
\item the same statement as (ii) with $\DD_1$ and $\DD_2$ interchanged.
\end{enumerate}
\end{lemma}
\begin{pf}
(i) is obvious from Lemma
\ref{_label_ifzcfacnwzeinxea}.
We prove (ii) and (iii).
By (\ref{_label_ftgqoznmugmufdat})
one can see that
(\ref{_label_mtknmljxgqmjihhu})
holds almost everywhere for fixed $f \in \DD_2$.
Let $A$ be a finite product of operators from
$\{(H-z)^{-1}\}_{z\in\RESOLVENT(H)}$.
By Lemma
\ref{_label_wbdmetvzzvzddjpo}
we see that
(\ref{_label_mtknmljxgqmjihhu})
remains valid when we replace
$f$ and $f_j$ by $Af$ and $Af_j$ respectively.
It follows that
one can make $\Lambda_0$ independent of the
choice of $f \in \DD_2$.
We can make $\Lambda_0$ smaller so that
(iii) holds as well.
This completes the proof.
\QED
\end{pf}

\begin{definition}
\label{_label_hpiiijxkslnbufyn}
Suppose that $H$ is a selfadjoint operator on a Hilbert space $\HilH$
and $\{\tau(\lambda)\}_{\lambda \in \RR}$ is a family
of pullback triples, where 
$\tau(\lambda) = (J(\lambda), F(\lambda), Y(\lambda))$
and the function spaces are as given by the commutative diagram below.
$$
\xymatrix{
&\DOM(Y(\lambda))
     \ar[ld]_{\iota}
     \ar[d]^{Y(\lambda)}
     \ar[rd]^{F(\lambda)Y(\lambda)}
\\
\HilH
&\ar[l]_{J(\lambda)} \XX(\lambda) \ar[r]^{F(\lambda)}
& \EE(\lambda)
}
$$
We call
$\{\tau(\lambda)\}_{\lambda \in \RR}$ a pullback representation
of the absolutely continuous part of $H$
if there exist a subset $\Lambda \subset \RR$
of full Lebesgue measure and
a span $\DD$
of an orthonormal basis by the resolvent of $H$
such that
\begin{enumerate}[i.]
\item $\DD$ is a core of $\tau(\lambda)$
for all $\lambda$;
\item for any $f \in \DD$ and $\lambda \in \Lambda$ we have
$$
\frac{d\|E_\lambda f\|^2}{d\lambda}
= \|F(\lambda)Y(\lambda)f\|_{\EE(\lambda)}^2;
$$
\item
for $\lambda \in \RR\setminus \Lambda$ we have
$\XX(\lambda) = \HilH$, $Y(\lambda)=J(\lambda)=1$,
 and $F(\lambda) = 0$.
\end{enumerate}
We say that two representations
$
\{\tau_j(\lambda)\}_{\lambda\in\RR}
, j=1,2
$
are equivalent if $\tau_1(\lambda)$ and $\tau_2(\lambda)$
are equivalent almost everywhere.
We refer to $\DD$ as a core of the representation.
\end{definition}

\begin{theorem}
Any selfadjoint operator on a separable Hilbert space
admits a pullback representation 
of the absolutely continuous part.  The representation is
unique up to equivalence.
\end{theorem}
\begin{pf}
Let $\DD$ and $\Lambda$ be as given by
Lemma \ref{_label_ifzcfacnwzeinxea}.
We see that
the function
$q_\lambda : \HilH \times \HilH \longrightarrow \CC$
defined for each $\lambda \in \Lambda$ by
\begin{equation}
\label{_label_xyqdnhcaipsrtncj}
q_\lambda(f, g) = \frac{d \la E_\lambda f | g \ra}{d \lambda},
\quad
f, g \in \DOM(q) = \DD
\end{equation}
is well-defined as a densely defined nonnegative quadratic
form on $\HilH$.  We set $q_\lambda = 0, \DOM(q_\lambda)=\DD$
for $\lambda \in \RR\setminus \Lambda$.
For each $\lambda \in \RR$ let
$$
J(\lambda) : \XX(\lambda) \longrightarrow \HilH,
\quad
\FF(\lambda) : \XX(\lambda) \longrightarrow \EE(\lambda),
\quad
Y(\lambda) : \HilH \longrightarrow \XX(\lambda)
$$
be a pullback triple of $q_\lambda$.
By Theorem
\ref{_label_dncrxfyquumwagzk}
we have
$$
q_\lambda(f, g) = \la F(\lambda) Y(\lambda)f | F(\lambda) Y(\lambda)g \ra_{\EE(\lambda)},
\quad f, g \in \XX(\lambda).
$$
We have thus proven that any selfadjoint operator admits
a pullback representation.
We now show the uniqueness.
Suppose that $\tau_j=\{\tau_j(\lambda)\}, j=1,2$
are pullback representations of $H$.
Write $\tau_j(\lambda) = (J_j(\lambda), F_j(\lambda), Y_j(\lambda))$.
In view of Lemma
\ref{_label_ekpkdvmkmvzchkga}
we may assume that 
$\DOM(Y_j(\lambda)) = \DD_j, j = 1,2$,
where each $\DD_j$ is a span of an orthonormal basis
by the resolvent of $H$.
Let $\Lambda_0$ be as in Lemma \ref{_label_brlgehxetoirlenm}.
Then for each $\lambda \in \Lambda_0$ there is
a well-defined linear operator
$Y_3(\lambda) : \HilH \longrightarrow \XX_1(\lambda)$ with
$\DOM(Y_3(\lambda))=\DD_2$ such that
\begin{equation}
\label{_label_wllyoawgiqwdtqwn}
J_1(\lambda)Y_3(\lambda)f = f,
\quad
\frac{
d\|E_\lambda f \|^2
}
{
d\lambda
}
= \| F_1(\lambda)Y_3(\lambda)f \|_{\EE_1(\lambda)}^2,
\end{equation}
\begin{equation}
\label{_label_oplfhldluqszxeis}
\frac{
d\la E_\lambda f | g \ra
}
{
d\lambda
}
= \la F_1(\lambda)Y_3(\lambda)f | F_1(\lambda) Y_1(\lambda)g \ra,
\quad
g \in \DD_1
\end{equation}
for all $f \in \DOM(Y_3(\lambda))$.
In view of Lemma \ref{_label_iucbdqsouvzpoyfk}
it suffices to show that $\tau_3(\lambda) = (J_1(\lambda), F_1(\lambda), Y_3(\lambda))$
is a pullback triple for all $\lambda \in \Lambda_0$.
Clearly $Y_3(\lambda)$ is densely defined and injective.
We now show that $Y_3(\lambda)$ has a dense range.
Let $\varphi \in \DD_1$.  Then there exists a 
sequence $\varphi_j \in \DD_2$ so that
\begin{equation}
\label{_label_lalnxsqcfuwfgjhf}
\|\varphi - \varphi_j\| \rightarrow 0,
\quad
\frac{
d\|E_\lambda (\varphi - \varphi_j)\|^2
}
{
d\lambda
}
\rightarrow 0
\end{equation}
for all $\lambda \in \Lambda_0$.
By
(\ref{_label_wllyoawgiqwdtqwn})
we see that $Y_3(\lambda)\varphi_j\rightarrow u$ in $\XX_1(\lambda)$
for some $u \in \XX_1(\lambda)$ and
(\ref{_label_oplfhldluqszxeis})
yields
$$
\la F_1(\lambda)[Y_3(\lambda)\varphi_j - Y_1(\lambda)\varphi]|
F_1(\lambda)Y_1(\lambda)g \ra
\rightarrow 0,
\quad
g \in \DD_1.
$$
Therefore $J_1(\lambda) u = \varphi$
and $F_1(\lambda) u = F_1(\lambda)Y_1(\lambda)\varphi$.
This means $u = Y_1(\lambda)\varphi$ and hence
$Y_3(\lambda)$ has a dense range.
We have thus shown that $Y_3(\lambda)$ is a completion operator.
This and 
(\ref{_label_wllyoawgiqwdtqwn})
means that $\tau_3(\lambda)$
is a pullback triple.
\QED
\end{pf}
  \section{Generalized eigenfunction expansion}
The purpose of this section is to show that the absolutely
continuous part of an arbitrary selfadjoint operator $H$
admits a generalized eigenfunction expansion.
We shall first
show that a pullback representation $\tau=\{\tau(\lambda)\}_{\lambda\in\RR}$
of $H$ induces a family of continuations
$\{H_\lambda\}_{\lambda \in \RR}$ of $H$.
Then we show that the spectral structue of
$H_\lambda$ and $H$ are identical except at energy $\lambda$
and the spectral density $E_\lambda'$ of $H$ at energy $\lambda$
is identical to a pullback of the orthogonal
projection onto the nullspace of $H_\lambda - \lambda$
for almost every $\lambda$.
This means that the structure of the
absolutely continuous part of
$H$ is encoded into the family
$\{\NULL(H_\lambda - \lambda)\}_{\lambda \in \RR}$
of nullspaces and we can recover the absolutely continous
part of the original operator from this family.

Throughout this section we fix a 
pullback representation
$\tau$
of the absolutely continuous part of $H$.
We use the same operators and function spaces as
in the Definition
\ref{_label_hpiiijxkslnbufyn}.

\begin{lemma}
\label{_label_ibkvqpiszhqhrxop}
Let $\lambda \in \RR$.  Then
\begin{enumerate}[i.]
\item
$f \in \DD$ and $\IM z \neq 0$ imply
\begin{equation}
\label{_label_whhfyoszafkklydb}
\FF(\lambda) Y(\lambda)(H-z)^{-1}f
= (\lambda-z)^{-1}\FF(\lambda) Y(\lambda)f;
\end{equation}
\item $(H-z)^{-1}\DD$ is a core of $\tau(\lambda)$ and we have
\begin{equation}
\label{_label_wtayvlxihsdneztb}
\|Y(\lambda)(H-z)^{-1}f \|_{\XX(\lambda)}
\leq |\IM z|^{-1} \|Y(\lambda)f\|_{\XX(\lambda)}
\end{equation}
for $f \in \DD$ and $\IM z \neq 0$;
\item $H$ is compatible with $Y(\lambda)$ in $\RESOLVENT(H)$.
\end{enumerate}
\end{lemma}
\begin{pf}
(i) is immediate from
Lemma
\ref{_label_wbdmetvzzvzddjpo}.
(ii).  From (i) we see that
$$
\|Y(\lambda)(H-z)^{-1}f\|^2_{\XX(\lambda)}
=\|(H-z)^{-1}f\|^2
+ |(\lambda-z)^{-1}|^2\|\FF(\lambda)Y(\lambda)f \|_{\EE(\lambda)}^2.
$$
This implies (\ref{_label_wtayvlxihsdneztb}).
Suppose that $u \in \XX(\lambda)$
and that
$
\la u | Y(\lambda)(H-z)^{-1}g \ra_{\XX(\lambda)} = 0
$
whenever $g \in \DD$.
It follows from
(i) that
$$
\la J(\lambda)u | J(\lambda)Y(\lambda)(H-z)^{-1}g \ra
+
(\lambda - \overline{z})^{-1}
\la F(\lambda)u | F(\lambda)Y(\lambda)g\ra
=0
$$
for $g \in \DD$.
Let $f_j$ be a sequence in $\DD$ such that $Y(\lambda)f_j \rightarrow u$
in $\XX(\lambda)$.  Note that $f_j \rightarrow J(\lambda)u$ in $\HilH$.
The inequality
(\ref{_label_wtayvlxihsdneztb})
implies
$$
\la f_j | (H-z)^{-1}f_j \ra
+
(\lambda - \overline{z})^{-1}
\| F(\lambda)Y(\lambda)f_j \|^2
\rightarrow 0
$$
and by taking the imaginary parts we see that
$$
\|f_j\|\rightarrow 0,
\quad
\|\FF(\lambda)Y(\lambda)f_j\|_{\EE(\lambda)}\rightarrow 0.
$$
Hence $u=0$.  By the Hahn-Banach theorem
we conclude that $Y(\lambda)(H-z)^{-1}\DD$ is dense in $\XX(\lambda)$.
(iii).
If $\lambda \in \RR \cap \RESOLVENT(H)$, then
$F(\lambda)=0$ and thus $Y(\lambda)$ is isometric.
This combined with (i) and (ii)
imply (i) through (iii) of Definition \ref{_label_zguvvrankpkeptbb}.
To verify the condition (iv)
suppose
$
Y(\lambda)f_j \rightarrow u
$
in $\XX(\lambda)$ with $f_j \in \DD$
and $Y(\lambda)(H-z)^{-1}f_j \rightarrow 0$ in $\XX(\lambda)$.
Then (i) yields
$$
\|(H-z)^{-1}f_j\|\rightarrow 0,\quad
\|\FF(\lambda)Y(\lambda)f_j\|_{\EE(\lambda)}\rightarrow 0
$$
so we must have $u = 0$.
This completes the proof.
\QED
\end{pf}

\begin{lemma}
\label{_label_zuhrnbgnefwrlnwp}
Let $\lambda \in \RR$.
Then $H$ has a selfadjoint continuation $H_\lambda$ to $\XX(\lambda)$
induced by $Y(\lambda)$
and it satisfies
\begin{equation}
\label{_label_fcgotxcqmnfmhrvi}
\FF(\lambda) \chi(H_\lambda)
= \chi(\lambda) \FF(\lambda),
\quad
J(\lambda)\chi(H_\lambda) = \chi(H) J(\lambda),
\end{equation}
for $\chi \in C_0^\infty(\RR)$.
\end{lemma}
\begin{pf}
Lemma \ref{_label_ibkvqpiszhqhrxop} implies that
$H$ has a continuation $H_\lambda$ to $\XX(\lambda)$
and that $\RESOLVENT(H) \subset \RESOLVENT(H_\lambda)$.
In view of Lemma
\ref{_label_ibkvqpiszhqhrxop}
the first identity of
(\ref{_label_fcgotxcqmnfmhrvi})
holds for $\chi(H) = (H-z)^{-1}$ with $z \in \RESOLVENT(H)$
and thus for all $\chi \in C_0^\infty(\RR)$.
Take $f_j \in \DD$ such that $Y(\lambda)f_j \rightarrow u$
in $\XX(\lambda)$.
Then we have
$$
J(\lambda)(H_\lambda - z)^{-1}Y(\lambda)f_j
= (H-z)^{-1}f_j
\rightarrow (H-z)^{-1}Ju
$$
for any $z \in \RESOLVENT(H)$.  This implies the second identity
of 
(\ref{_label_fcgotxcqmnfmhrvi}).
For any $z \in \RESOLVENT(H)$ we have
$$
\IM \la (H_\lambda-z)^{-1} u | u\ra_{\XX(\lambda)}
= \IM z \| (H_\lambda-z)^{-1}u\|_{\XX(\lambda)}^2
$$
for all $u \in Y(\lambda)\DD$ and thus for all $u \in \XX(\lambda)$.
Therefore $H_\lambda$ is symmetric and hence
selfadjoint by Corollary
\ref{_label_qahdlsdgvfgxtoqj}.
\QED
\end{pf}

By Lemma \ref{_label_zuhrnbgnefwrlnwp}
$J(\lambda)\NULL(F(\lambda))$ is dense in $\HilH$
for all $\lambda \in \RR$, so we have the following
\begin{corollary}
The pullback triple $\tau(\lambda)$
is splitting for all $\lambda \in \RR$.
\end{corollary}

\begin{theorem}
\label{_label_smvnkbchdupyvpxz}
For any $\lambda \in \RR$
the spectrum of $H$ and $H_\lambda$ are identical.
We have
\begin{equation}
\label{_label_fgepeuhrmmmqocya}
J(\lambda) H_\lambda
= HJ(\lambda),
\quad
\FF(\lambda) H_\lambda
\subset \lambda \FF(\lambda).
\end{equation}
$J(\lambda)$ and $\FF(\lambda)$
are partial isometries.
For $\lambda \in \Lambda$ we have
\begin{gather}
\label{_label_mfnwnznlastdxzmn}
\CLOSURE{\RANGE(H_\lambda - \lambda)}
= \NULL(\FF(\lambda))
= \RANGE(J(\lambda)^*),
\quad
\RANGE(J(\lambda)) = \HilH,
\\
\label{_label_ljtvvbiirzueboex}
\NULL(H_\lambda - \lambda)
= \NULL(J(\lambda))
= \RANGE(F(\lambda)^*),
\quad
\RANGE(F(\lambda)) = \EE(\lambda).
\end{gather}
\end{theorem}
\begin{pf}
Lemma \ref{_label_zuhrnbgnefwrlnwp} implies
(\ref{_label_fgepeuhrmmmqocya}).
We now prove
\begin{equation}
\label{_label_zofbebgqyddxiqgv}
\NULL(H_\lambda-\lambda)
= \NULL(J(\lambda)).
\end{equation}
As $H_\lambda$ is selfadjoint
and $\FF(\lambda)(H_\lambda - \lambda)v = 0$ for any $v \in \DOM(H_\lambda)$,
$u \in \NULL(H_\lambda - \lambda)$ if and only if
$
\la J(\lambda)u | J(\lambda)(H_\lambda - \lambda) v \ra = 0
$
for all $v \in \DOM(H_\lambda)$.
This condition is equivalent to $J(\lambda)u = 0$
because $J(\lambda)(H_\lambda-\lambda)$ has a dense range.
Hence we obtain (\ref{_label_zofbebgqyddxiqgv}).
By Proposition \ref{_label_cqwlrlcwckueubdr},
(\ref{_label_mfnwnznlastdxzmn})
and
(\ref{_label_ljtvvbiirzueboex})
follow.
We now prove $\SPEC(H) = \SPEC(H_\lambda)$.
This is obvious if $\lambda \in \RR\setminus\Lambda$,
so we may assume $\lambda \in \Lambda$.
Then
$J(\lambda)$ is a partial isometry
with initial space $\NULL(Q(\lambda))$ and final space $\HilH$,
where $Q(\lambda)=\FF(\lambda)^*\FF(\lambda)$.
Note that $Q(\lambda)$
is the orthogonal projection onto
the eigenspace of $H_\lambda$ with eigenvalue $\lambda$.
This and 
(\ref{_label_fgepeuhrmmmqocya})
imply that $E_H(I)$ is unitarily equivalent to $E_{H_\lambda}(I\setminus\{\lambda\})$
for any interval $I$.
If $\lambda \in \RESOLVENT(H)$,
then $E_{H_\lambda}(\{\lambda\}) = Q(\lambda) = 0$ and thus
$\lambda \in \RESOLVENT(H_\lambda)$.
Conversely, $\lambda \in \RESOLVENT(H_\lambda)$ implies
$Q(\lambda) = 0$ and so $\lambda \in \RESOLVENT(H)$.
This completes the proof.
\QED
\end{pf}

\begin{theorem}
\label{_label_bzlmyicsrcbuibbj}
There exists a partial isometry
\begin{equation}
\label{_label_pplkqyrlkcbdkxfd}
\FF : \HilH \longrightarrow \int_{\RR}^\oplus \EE(\lambda) d\lambda
\end{equation}
such that
$$
(\FF f)(\lambda) = \FF(\lambda)Y(\lambda)f,
\quad
f \in \DD,
\quad
\lambda \in \RR.
$$
We have
\begin{equation}
\label{_label_mdwssftxbxynfpcn}
F^*F = P_\AC(H),
\quad
FF^* = 1.
\end{equation}
\end{theorem}
\begin{remark}
We remark that
$$
\frac{
d\la E_\lambda f | g \ra
}{
d\lambda
}
= \la Q(\lambda)Y(\lambda)f | Y(\lambda)g\ra_{\XX(\lambda)},
\quad
f, g \in \DD
$$
almost everywhere, where $Q(\lambda) = \FF(\lambda)^*\FF(\lambda)$.
It follows that the spectral density $E_\lambda'$
is a pullback of
the eigenprojection
$Q(\lambda)$ of $H_\lambda$ by the completion operator $Y(\lambda)$ and Theorem
\ref{_label_bzlmyicsrcbuibbj}
implies
$$
\la H_{\AC} f | g \ra
= \int_{\RR}
\lambda
\la Q(\lambda)Y(\lambda)f | Y(\lambda)g\ra_{\XX(\lambda)}d\lambda.
\quad
$$
The absolutely continuous part of 
a selfadjoint operator thus admits
a generalized eigenfunction expansion.
It should be noted that $Q(\lambda)=0$
if and only if the spectral density
is closable at $\lambda$ and
that the size of the generalized eigenspace 
$$
\NULL(H_\lambda - \lambda)
= \RANGE(Q(\lambda))
= \NULL(J(\lambda))
$$
is determined by how far is
the spectral density from being closable
at energy $\lambda$.
\end{remark}
\begin{pf}
It is easy to see that
$$
\|P_{\AC}(H) f\|^2
=\int_\RR \|\FF(\lambda)Y(\lambda)f\|_{\EE(\lambda)}^2 d\lambda
\quad
f \in \DD.
$$
It follows that the partial isometry $F$ exists and we have
$F^*F = P_\AC(H)$.
We now show that $\FF$ is surjective.
By the Hahn-Banach theorem it suffices to prove that
\begin{gather*}
\int_\RR
\la 
g(\lambda)
|
\FF(\lambda)Y(\lambda)f
\ra_{\EE(\lambda)}
d\lambda = 0,
\quad
f \in \DD,
\\
\int_\RR \|g(\lambda)\|_{\EE(\lambda)}^2 d\lambda < \infty
\end{gather*}
implies $g(\lambda) = 0$ for almost every $\lambda$.
The identity
(\ref{_label_fcgotxcqmnfmhrvi})
yields
$$
\int_\RR \varphi(\lambda)
\la g(\lambda) | \FF(\lambda)Y(\lambda)f \ra_{\EE(\lambda)}
d\lambda = 0
$$
for any $\varphi \in C_0^\infty(\RR)$ and $f \in \DD$.
Since $\DD$ is the span of an orthonormal basis by the
resolvents of $H$, there exists a set $\Lambda_0\subset\RR$
of full Lebesgue measure which is independent of $f$
so that
$\la g(\lambda) | \FF(\lambda) Y(\lambda)f\ra_{\EE(\lambda)} = 0$
for $\lambda \in \Lambda_0$ and $f \in \DD$.
Since $\FF(\lambda)$ is bounded and surjective,
we must have $g(\lambda) = 0$ for almost every $\lambda$.
This completes the proof.
\QED
\end{pf}
  \section{Acknowledgements}
The author is grateful to Professor Shu Nakamura
for valuable discussions and constructive remarks.





\bibliographystyle{model1b-num-names}
\bibliography{Main}







\end{document}